\documentclass[11pt]{article}
\usepackage{amsmath,amssymb,bbm,color}
\font\nineam=msam10
\newcommand{\pdwJump}{{\!\hbox{\nineam\char'110}}}

\newcommand{\nobracket}{}\newcommand{\nocomma}{} \newcommand{\nosymbol}{}   \newcommand{\tmem}[1]{{\em #1\/}} \newcommand{\tmmathbf}[1]{\ensuremath{\boldsymbol{#1}}} \newcommand{\tmop}[1]{\ensuremath{\operatorname{#1}}} \newcommand{\tmrsup}[1]{\textsuperscript{#1}} \newcommand{\tmsamp}[1]{\textsf{#1}} \newcommand{\tmstrong}[1]{\textbf{#1}}  \newcommand{\tmtextbf}[1]{{\bfseries{#1}}} \newcommand{\tmtextit}[1]{{\itshape{#1}}} \definecolor{grey}{rgb}{0.75,0.75,0.75} \definecolor{orange}{rgb}{1.0,0.5,0.5} \definecolor{brown}{rgb}{0.5,0.25,0.0} \definecolor{pink}{rgb}{1.0,0.5,0.5} \newtheorem{definition}{Definition} \newtheorem{lemma}{Lemma} \newtheorem{theorem}{Theorem} 
\newtheorem{corollary}{Corollary}

\def\re{\mathbbm{R}}
\def\da{\downarrow\!}

\def\qed{Q.E.D.}
\def\bm{\boldmath}

\title{  \large {Discrete Transfinite Computation} }

\author{ \normalsize P. D. Welch\\
 \small School of Mathematics, University of Bristol, England. } 
\date{ \small }%\today}
\begin{document}

\maketitle

\abstract{We describe various computational models based initially, but not exclusively, on that of the Turing machine, that are generalized to allow for transfinitely many computational steps.  Variants of such machines are considered that have longer tapes than the standard model, or that work on ordinals rather than numbers. We outline the connections between such models and the older theories of recursion in higher types, generalized recursion theory, and recursion on ordinals such as $\alpha$-recursion. We conclude that, in particular, polynomial time computation on $\omega$-strings is well modelled by several convergent conceptions.}

\section{Introduction}

There has been a resurgence of interest in models of infinitary computation in
the last decade. I say resurgence because there has been for 50 years or more
models of computation that deal with sets of integers, or objects of even
higher type, in generalized recursion theory.
%which must be recognised as infinitary in some form. 
Such a theory was initiated by Kleene's generalization
of his equational calculus for ordinary recursion theory {\cite{Kl59}}-{\cite{Kl63}}. Whilst that was indeed  a generalized recursion theory
some commentators remarked that it was possible to view, for example,
what came to be called Kleene recursion, as having a machine-like model. The
difference here was that the machine would have a countable memory, a
countable tape or tapes, but an ability to manipulate that countable memory in
finite time, or  equivalently  in one step. Thus a query about a set of
integers $r$ say, say coded as an element of Cantor space, as to whether $r$ was an element of an oracle set $A \subseteq 2^{\mathbbm{N}}$, could be
considered a single computational operation. Moreover  $r$ could be moved from
a tape to a storage or otherwise manipulated or altered in infinitely many
places at once by one step of some procedure. This is what Addison and Spector
called the ``$\aleph_0$-mind''  (see{\cite{Rog}}, p.405;  further Kreisel spoke of a generalized Church's thesis). We thus had essentially a
generalization of the original Turing machine model. It should be noted
however (in contradistinction to the transfinite models of more recent years)
that the computations would all be represented by wellfounded computation
trees when convergent, and would otherwise be deemed non-convergent: an
infinitely long linear sequence of computations represents failure. By declaring that countably
many operations on a set of integers could be done in ``finite time'' or all at
once and so count as one step, one simply sidesteps the difficulty of thinking
of infinitely many steps of a successful computation that is yet to be
continued.

Thus in the past mathematicians have been more inclined to consider
wellfounded computation but applied to infinite objects, rather than
considering transfinite computations containing paths of transfinite
{\tmem{order type}}. Examples of the former: Kleene recursion as
HYP-recursion, or the Blum-Shub-Smale (BSS, {\cite{BlShSm89}}) machines acting
with real inputs (or indeed other mathematical objects taken from a ring with
suitable functions).  One could almost say that a hallmark of the
generalizations of recursion theory to ``higher'', or
``generalized'' recursion theory has been that it has considered {\tmem{only}}
wellfounded computation albeit with infinitary objects of higher
types. Sacks's $E$-recursion {\cite{Sa90}} again considers essentially the same
paradigm of building up functions applied now to sets in general: a
non-convergent computation is exemplified by an infinite path in a computation
tree, thus rendering the tree ill-founded.

However what if we were nevertheless to think also of transfinite stages of computation?
There is a natural reaction to this: we feel that ``computation'' is so tied
up with our notions of finitary algorithm, indeed effective algorithm, that it
is too heterodox to consider transfinitely numbered steps as `computation.'
However why should this be? \ We are probably constrained to think of well
defined computations as shying away from anything that smells of a
{\tmem{supertask}}: the quandaries arising from Thomson's Lamp {\cite{Th54}}
seem to have used up a surprising lot of space and  ink, for
what is after all a simple definitional problem. So supertasks have been
banished from computational study.

{\tmsamp{}}However transfinite recursion or procedures are not at all alien,
if not common in the everyday experience of the mathematician - there are
after all, few $\Pi^1_1$-complete sets occuring as `natural' mathematical
objects - but the early canonical example arises in the Cantor-Bendixson
process. Here with a simple decreasing chain of sets 
$X_{\alpha} \supseteq X_{\alpha + 1}$ and intersections at limits:
$X_{\lambda} = \cap_{\alpha < \lambda} X_{\alpha}$, one has a monotone process.  Again the monotonicity phenomenon occurs centrally in Kleene's recursion
theory of higher types, and we feel safer with a monotone process than a
non-monotone or discontinuous one.

Notwithstanding these qualms, the current chapter reviews the recent
descriptions of various machine models, including that of Turing's original
machine itself, which can be given a defined behaviour at limit stages of
time, enabling them to compute through recursive ordinals and beyond. This behaviour, apart from a few very elementary models, is signified by being non-monotonic, or {\em quasi-inductive.}

We shall see that the various models link into several areas of  modern
logic: besides recursion theory, set theory and the study of subsystems of second order analysis play a role. Questions arise concerning the strengths of models that operate at the level one type above that
of the integers. This
may be one of ordinal types: how long a well ordered sequence of steps must a
machine undertake in order to deliver its output? Or it may be of possible output:
if a machine produces real numbers, which ordinals can be coded as output
reals? And so on and so forth.

\subsection{Computation in the limit}

To start at the beginning, steps towards transfinite considerations, or at
least that of considering what might have occurred on a tape at a finite stage,
come immediately after considering the halting problem. The universal Turing machine
can be designed to print out on an infinite output tape the code numbers $e$
of programmes $P_e$ that will halt on input $e$: thus $P_e (e)\! \da
\nosymbol$ . This $\Sigma^0_1$ set, as is well known, is {\tmem{complete}}:
any other $\Sigma^0_1$ is (ordinary) Turing reducible to it.

Putnam {\cite{Pu65}} (and Gold {\cite{Go65}}) went a step further.

\begin{definition}
  {\tmem{({\cite{Pu65}})}} $P$ is a {\tmem{trial and error predicate}}
  if and only if there is a general recursive function $f$ such that for every
  $x_1, \ldots, x_n$:\\

 $
  \begin{array}{rcl}
  P (x_1, \ldots, x_n) & \equiv &\lim_{y \rightarrow \infty} f ( x_1,
  \ldots, x_n, y) = 1 \\
  \neg P ( x_1, \ldots, x_n)& \equiv & \lim_{y \rightarrow \infty} f (
  x_1, \ldots, x_n, y) = 0.
\end{array}
$
\end{definition}

Running such a precedure on a Turing machine allows us to print out a
$\Delta^0_2$ set $A$'s characteristic function on the output tape. In order to
do this we are forced to allow the machine to change its mind about $n \in A$
and so repeatedly substitute a 0 for a 1 or vice versa in the $n$'th cell of
the output tape. However, and this is the point, at most finitely many changes
are to be made to that particular cell's value. 
It is this feature of not knowing at any finite given time whether further alterations are to made, that makes this a transition from a computable set to a non-computable one.

By a recursive division of the
working area up into infinitely many infinite pieces, one can arrange for the
correct computation of all ?$m \in A ?$ to be done on the one machine, and the
correct values placed on the output tape.

However this is as far as one can go if one imposes the (very obvious,
practical) rule that a cell's value can only be altered finitely often. In
order to get a $\Sigma^0_2$ set's characteristic function written to the output
tape, then in general one cannot guarantee that a cell's value is changed
finitely often. Then immediately one is in the hazardous arena of supertasks.

Nevertheless let us play the mathematicians' game of generalizing for generalization's
sake: let us by fiat declare a cell's value that has switched infinitely often
$0 \rightarrow 1 \rightarrow 0$ to be $0$ at ``time $\omega$''. \ With this
$\liminf$ declaration one has, mathematically at least, written down the
$\Sigma^0_2$-set on the output tape, again at time $\omega$.

Following this through we may contemplate continuing the computation at times
$\omega + 1, \omega + 2, \ldots, \omega + \omega \nocomma \nocomma \ldots$. \
Let $\langle C_i (\alpha) \mid i \in \omega, \alpha \in \tmop{On} \rangle$
denote the contents of cell $C_i$ at time $\alpha$ in $\tmop{On}$. Let $l
(\alpha) \in \omega$ represent the cell number being observed at time
$\alpha$. Similarly let $q (\alpha)$ denote the state of the machine/program
at time $\alpha$.  We merely need to specify (i) the read/write head position $l
(\lambda)$ and (ii) the next state $q (\lambda)$ for limit times $\lambda$
(whilst $l (\alpha + 1)$ and $q (\alpha + 1)$ are obtained by just following
the usual rules for head movement and state change according to the standard Turing
transition table. (Note that the stock of programmes has not changed, hence $q
(\alpha)$ names one of the finitely many states of the usual transition table;
we are merely enlarging the possible behaviour.)

\def\col#1#2#3{\hbox{\vbox{\baselineskip=0pt\parskip=0pt\cell#1\cell#2\cell#3}}}

\newcommand{\cell}[1]{\boxit{\hbox to 17pt{\strut\hfil$#1$\hfil}}}

\newcommand{\head}[2]{\lower2pt\vbox{\hbox{\strut\footnotesize\it\hskip3pt#2}\boxit{\cell#1}}}

\newcommand{\boxit}[1]{\setbox4=\hbox{\kern2pt#1\kern2pt}\hbox{\vrule\vbox{\hrule\kern2pt\box4\kern2pt\hrule}\vrule}}

\newcommand{\Col}[3]{\hbox{\vbox{\baselineskip=0pt\parskip=0pt\cell#1\cell#2\cell#3}}}

\newcommand{\tapenames}{\raise 10pt\vbox to 40pt{\hbox to .8in{\it\hfill Input: \strut}\vfill\hbox to
.8in{\it\hfill Scratch: \strut}\vfill\hbox to .8in{\it\hfill Output: \strut}}}

\newcommand{\Head}[4]{\lower2pt\vbox{\hbox to25pt{\strut\footnotesize\it\hfill#4\hfill}\boxit{\Col#1#2#3}}}

\newcommand{\Dots}{\raise 7pt\vbox to 50pt{\hbox{\ $\cdots$\strut}\vfill\hbox{\ $\cdots$\strut}\vfill\hbox{\ $\cdots$\strut}}}

\begin{figure}[h]

$$\tapenames\Col101\Col110\Head000{R/W}\Col100\Col111\Col011\Col000\Col001\Col000\Dots$$
{\caption{\label{Fig4} A 3-tape Infinite Time Turing Machine}}
\end{figure}

We intend that control of the program at a limit time $\lambda$, be placed at
the beginning of the outermost loop, or subroutine call, that was called
unboundedly often below $\lambda$. We thus set:

$$
q (\lambda)   =  \tmop{Liminf}_{\alpha <\lambda} q(\alpha).$$

For limit $\lambda$ we set $C_i (\lambda) $ by:\\

$\begin{array}{rcl}C_i (\lambda) &=&  k \mbox{ if }\exists \alpha < \lambda \forall
\beta < \lambda (\alpha < \beta \longrightarrow C_i (\beta) = k)\mbox{ for } k
\in \{0, 1\};\\

 & = & 0 \mbox{ otherwise.}
\end{array}
$\\

The R/W head we place according to the above, also using a modified Liminf
rule:\\

$
\begin{array}{lcl} 
l (\lambda) & = & \tmop{Liminf}^\ast\langle l(\alpha)\mid \alpha <\lambda\rangle
\end{array}
$\\

This is not exactly the arrangement that Hamkins and Lewis specified in
{\cite{HL}} but it is inessentially different from it (HL specified a
special {\tmem{limit state}} $q_{\lambda}$ which the machine entered into
automatically at limit stages, and the head was always set back to the start
of the tape. They specified (which we shall keep here) that the machine be a
three tape machine.

Input then can consist of a set of integers, suitably coded as an element of
2\tmrsup{$\mathbbm{N}$} \ on $\langle C_{3 i} \rangle_i$ and output likewise
is such an element on $\langle C_{3 i + 2} \rangle_i$. Thus there is little
difference in a machine with an oracle $Z \subseteq \omega$ and one acting on
input $Z$ coded onto the input tape. \ However we immediately see the
possibility of higher type computation: we may have some $\mathcal{Z}
\subseteq 2^{\mathbbm{N}}$ and then we add a query state which asks if, say,
the scratch tape $\langle C_{3 i + 1} \rangle_i$'s contents is or is not an
element of $\mathcal{Z}$.

We have thus completely specified the ITTM's behaviour.  The scene is thus
set to ask what such machines are capable of.  We defer discussion of this until
Section 2, whilst we outline the rest of this chapter here.

In one sense we have here a logician's plaything: the Turing model has been
taken over and redesigned with a heavy-handed liminf rule of behaviour. This
liminf operation at limit stages is almost tantamount to an infinitary logical
rule, and most of the behaviour the machine exhibits is traceable to this
rule.  But then of course it  has to be, what else is there? Nevertheless
this model and those that have been studied subsequently have a number of
connections or aspects with other areas of logic. Firstly, with weak
subsystems of analysis: it is immediately clear that the behaviour of such
machines is dependent on what ordinals are available. A machine may now halt
at some transfinite stage, or may enter an infinitely repeating loop; but any
theory that seeks to describe such machines fully is a theory which implies
the existence of sufficiently long wellorderings along which such a machine
can run (or be simulated as running). We may thus ask ``What (sub-) system of
analysis is needed in which to discuss such a machine''? \ We shall see that
machine variants may require longer or shorter wellorderings, thus their 
theory can be discussed within different subsystems.

Secondly, we can ask how the computable functions/sets of such a model fit in
with the earlier theories of generalized recursion theory of the 1960's and
70's.  For example there is naturally associated with ITTM's a so-called
{\tmem{Spector Class}} of sets. Such classes arise canonically in the
generalized recursion theories of that era through notions of definability.

Once one model has been defined it is very tempting to define variants. One
such is the {\tmem{Infinite Time Register Machine}} (ITRM's - due to Koepke
{\cite{KoMi08}}) which essentially does for Shepherdson-Sturgis machines what
HL does for Turing machines. \ Whilst at the finite level these two models are
equal in power, their infinitary versions differ considerably, the ITTM's
being much stronger.  The ITRM model is discussed in Section 3. 

Just as for ordinary recursion on $\omega$ the TM model with a putative tape
in order type $\omega$ length is used, so when considering notions of
$\alpha$-recursion theory for admissible ordinals $\alpha$, it is possible to
think of tapes also unfettered by having finite initial segments: we may
consider machines with tapes of order type $\alpha$ and think of computing
along $\alpha$ witb such machines. What is the relation to this kind of
computation and $\alpha$-recursion theory?

One can contemplate even machines with an $\tmop{On}$-length tape. \ It turns
out (Koepke {\cite{KK06}}) that this delivers a rather nice presentation of
G{\"o}del's constructible hierarchy. Finally discussed here is the notion of a
{\tmem{Blum-Shub-Smale}} machine ({\cite{BlShSm89}} acting transfinitely. With
some continuity requirement imposed on register contents for limit times, we
see that functions such as exponentiation $e^x$ which are not BSS computable,
become naturally IBBS computable. Moreover there is a nice equivalence between
their decidable reals, and those produced by the {\tmem{Safe Set Recursion} (``SSR'')}
of Beckmann, Buss, and S. Friedman, which can be thought of as generalizing to
transfinite sets notions of polynomial time computable functions on integers.
\ Briefly put, a polynomial time algorithm using $\omega$ as an input string,
should be halting by some time $\omega^n$ \ for some finite $n$. The IBBS
computable reals are then identical to the SSR-computable reals. The
background second order theory needed to run IBBS machines lies intermediate
between $\tmop{WKL}_0 $ and $\tmop{ATR}_0$.

The relation of ITTM's to Kleene recursion is discussed in Section 2.

\section{\tmstrong{What ITTM's can achieve.}}

Hamkins and Lewis in {\cite{HL}} explore at length the properties of ITTM's:
they demonstrate the natural notion of a {\tmem{universal}} such machine, and
hence an $S^n_m$-theorem and the Recursion Theorems. A number of questions
immediately spring to mind:\\

{\tmem{Q.}} \ What is the halting set $H = \{ e \in \omega \mid P_e (0)\!\da \,\}$?\\

Here $\langle P_e \rangle_e$ enumerates the usual Turing machine
programs/transition tables (and we use \ $P_e (x) \da y$ to denote that
the $e$'th program on input $x \in \mathbbm{N}$ or in $2^{\mathbbm{N}}$ halts
with output $y$. (If we are unconcerned about the $y$ we omit reference to
it.) An ITTM computation such as this can now halt in $\omega$ or more many
steps. But how long should we wait to see if $P_e (0)\da$ or not? \
This is behind the following definitions.

\begin{definition}
  (i) We write `$P_e (n) \!\da^{\alpha} y$' if $\text{$P_e (n) \!\da
  y$}$ in exactly $\alpha$ steps. We call $\alpha$ {\tmem{clockable}} if
  $\exists e \exists n\in
  \omega \exists y$ $P_e (n) \!\da^{\alpha} y$.

  (ii) A real $y \in 2^{\mathbbm{N}}$ is {\tmem{writable}} if there are $n, e\in \omega$
  with $\text{$P_e (n) \!\da y$}$; an ordinal $\beta$ is called writable,
  if $\beta$ has a writable code $y$.
\end{definition}

We may consider a triple
$s(\alpha)= \langle l (\alpha), q (\alpha), \langle C_i (\alpha)_i \rangle \rangle$ 
as a  {\tmem{snapshot}} of a machine at time $\alpha$, 
which contains all the relevant information at that moment. A {\em computation} is
then given by a wellordered sequence of snapshots. There are two possible
outcomes: there is some time $\alpha$ at which the computation{\tmem{ halts}},
or else there must be some stage $\alpha_0$ at which the computation enters
the beginning of a loop, and from then on throughout the ordinals it must
iterate through this loop. It is easy either by elementary arguments or simply
by L\"owenheim-Skolem, to see that such an $\alpha_0$ must be a countable ordinal, and
moreover that the periodicity of the cycling loop is likewise countable.

The property of
being a `{\tmem{well-ordered sequence of snapshots in the computation $P_e
(x)$'}} is $\Pi^1_1$ as a relation of $e$ and $x$.  Hence `$P_e (x) \!\da
y$' is $\Delta^1_2$:\\

$\exists w (w$ {\tmem{codes a halting computation of $P_e (x)$, with $y$
written on the output tape at the final stage$)$}} $\Longleftrightarrow$

$\forall w (w$ {\tmem{codes a computation of $P_e (x)$ that is either halting
or performs a repeating infinite loop $\longrightarrow w$ codes a halting
computation with $y$ on the output tape}}$.)$\\

Likewise $P_e (x)\! \uparrow$ is also $\Delta^1_2$.
By the above discussion then it is immediate that the clockable and writable
ordinals are all countable. Let $\lambda =_{\tmop{df}} \sup \{ \alpha \mid
\alpha \nobracket$ is writable$\} \nobracket$; let $\gamma =_{\tmop{df}} \sup
\{ \alpha \mid \alpha \nobracket$ is clockable$\} \nobracket$. Hamkins-Lewis
showed that $\lambda \leq \gamma$.\\

{\tmem{Q2}} Is $\lambda = \gamma$?\\

\begin{definition}
  \label{haltings} (i) $x^{\nabla} =\{e | P_e (x)
  \!\downarrow\}$ (The halting set on integers).

  (ii) $X^{\pdwJump} =^{} \{(e, y) | P_e^X (y) \!\downarrow\}$ (The halting
  set on reals relativised to $X \subseteq 2^{\mathbbm{N}}$).
\end{definition}

This yields the halting sets, both for computations on integers and secondly on reals where by the latter we include the instruction for the ITTM to query whether the current scratch tape's contents considered as a real, is in $X$.

\begin{definition}
  (i) $R (x)$ is an $\tmop{ITTM}$-{\tmem{semi-decidable}} predicate if there
  is an index $e$ so that:
$$\forall
  x (R (x) \leftrightarrow P_e (x) \!\da 1)$$
 \noindent (ii) A predicate $R$ is $\tmop{ITTM}${\tmem{-decidable}} if both $R$ and
  $\neg R$ are $\tmop{ITTM}$-semi-decidable.
\end{definition}

{\tmem{Q3 }} What are the ITTM-(semi-)decidable sets of integers, or reals? What is $x^\nabla$?\\

The last question seems somewhat impenetrable without a characterisation of
the halting behaviour of ITTM's - and one version of that problem is Question 2.
To analyse decidability one route might be through a version of Kleene's
{\tmem{Normal Form Theorem}} in the context of ITTM's. However there is an
obvious type difference: a successful halting computation of an ordinary
Turing machine in a finite amount of time can be coded as a finite sequence of
finite ordinary-TM-snapshots, and thus through the usual coding devices, by an
integer. This is essentially the heart of Kleene's $T$-predicate argument.
Thus for standard TM's Kleene demonstrated that given an index $e$ we may
effectively find an index $e'$ so that for any $n$: \ $P_e (n) \da \,\,\longrightarrow P_{e'} (n) \da M$ \ where $M \in \mathbbm{N}$ is a code
for the whole computation sequence of $P_e (n)$.

However this is too simple here: $P_e (0)\!\da$ may halt at some
transfinite time $\beta \geq \omega$. Hence the halting computation is only
codable by an infinite sequence of infinite ITTM-snapshots $S$ say, of some
order type $\tau$. In order for there to be a chance of having another index
$e'$ with $P_{e'} (n) \da y$ where $y$ codes such a sequence $S$, one
has to know at the very least that there is an $e_0$ so that $P_{e_0} (n)$
halts with output a code for $\tau$. In other words that $\tau$ be
{\tmem{writable}}. \ \ Thus we need an affirmative answer to {\tmem{Q2}} at
the very least.

Interestingly the key to answering halting behaviour is not to aim straight
for an analysis of halting {\tmem{per se}}, but at another phenomenon that is
peculiarly significant to ITTM's. There can be computations $P_e (x)$ that
whilst they have not formally halted, nevertheless from some point in time
onwards, leave their output tapes alone, and just churn around for ever
perhaps doodling or making entries on their scratch tape. The output has
however {\tmem{stabilized}}. \ We formally define this as follows:

\begin{definition}
  (i) Suppose for the computation $P_e (x)$ the machine does not halt then we
 write $P_e (x)\! \uparrow$ ; if eventually the output tape does have a stable value
  $y \in 2^{\mathbbm{N}}$ then we write: $P_e (x)\! \uparrow y$ and we say that
  $y$ is {\tmem{eventually $x$-writable}}.

  (ii) $R (x)$ is an {\tmem{eventually
  }}$\tmop{ITTM}$-{\tmem{semi-decidable}} predicate  if there
  is an index $e$ so that:
 $$\forall
  x (R (x) \leftrightarrow P_e (x) \!\uparrow 1)$$
  (iii) A predicate $R$ is {\tmem{eventually
  }}$\tmop{ITTM}${\tmem{-decidable}} if both $R$ and $\neg R$ are eventually
  $\tmop{ITTM}${\tmem{-{\tmem{semi-decidable.}}}}
\end{definition}

This proliferation of notions is not gratuitous: it turns out that answering
{\tmem{Q2}} on clockables {\tmem{vis {\`a} vis}} writables, depends on
successfully analysing stabilization patterns of individual cells $C_i$ during
the course of the computation $P_e (n)$. \ The moral is that stabilization is
anterior to halting. The following lemma illustrates the point.

\begin{theorem}
  \label{lzs-II}{\tmem{{\tmstrong{(The {\bm$\lambda, \zeta,
  \Sigma$}-Theorem)}} (Welch {cf. \cite{W}, \cite{W09}})}} (i) Any ITTM computation $P_e (x)$ which
  halts, does so by time $\lambda^x$, the latter being defined as the supremum of
  the $x$-writable ordinals.

  {\noindent}(ii) Any computation $P_e (x)$ with eventually stable output
  tape, will stabilize before the time $\zeta^x$ defined as the supremum of the
  eventually x-writable ordinals.

  {\noindent}(iii) Moreover $\zeta^x$ is the least ordinal so that there
  exists $\Sigma^x > \zeta^x$ with the property that
  \[ \text{$L_{\zeta^x} [x] \prec_{\Sigma_2} L_{\Sigma^x} [x] ;$} \]
  (iv) Then $\lambda^x$ is also characterised as the least ordinal satisfying:
  \[ \text{$L_{\lambda^x} [x] \prec_{\Sigma_1} L_{\zeta^x} [x] .$} \]
\end{theorem}

If we unpack the contents here, answers to our questions are given by (iii)
and (iv). Let us take $x = \varnothing$ \ so that we may consider the
unrelativised case. Our machine-theoretic structure and operations are highly
absolute and it is clear that running the machine inside the constructible
hierarchy of $L_{\alpha}$'s yields the same snapshot sequence as considering
running the machine in $V$. If $P_e (n) \da$ then this is a
$\Sigma_1$-statement (in the language of set theory). As halting is merely a
very special case of stabilization, then we have that
$$P_e
(n) \da \,\,\leftrightarrow (P_e (n)\da\,\,)^{L_{\zeta}}
\leftrightarrow (P_e (n) \da\,\,)^{L_{\lambda}}$$

\noindent (the latter because $L_{\lambda} \prec _{\Sigma_1} L_{\zeta}$). \ Hence the
computation must halt before $\lambda$. \ Hence the answer to {\tmem{Q2}} is
affirmative: every halting time (of an integer computation) is a writable
ordinal. \ One quickly sees that a set of integers is ITTM-decidable if and
only if it is an element of $L_{\lambda}$. \ It is ITTM-semi-decidable if and
only if it is $\Sigma_1 (L_{\lambda})$.

Since the limit rules for ITTM's are intrinsically of a $\Sigma_{2}$-nature, with hindsight it is perhaps not surprising that this would feature in the $(\zeta,\Sigma)$ pair arising as they do: after all the snapshot of the universal ITTM at time $\zeta$ is going to be coded into the $\Sigma_{2}$-Theory of this $L_{\zeta}$.  The universality of the machine is then apparent in the fact that by stage $\zeta$ it will have ``constructed'' all the constructible sets in $L_{\zeta}$.\\

As a corollary one obtains then:

\begin{theorem}
  \label{Normal2}{\tmem{\tmtextbf{(Normal Form Theorem)}(Welch)}} (a) For any ITTM
  computable function $\varphi_e$ we can effectively find another ITTM
  computable function $\varphi_{e'}$ so that on any input $x$ from
  $2^{\mathbbm{N}}$ if $\varphi_e (x) \!\downarrow$ then $\varphi_{e'} (x)
  \!\downarrow y \in 2^{\mathbbm{N}}$, where $y$ codes a wellordered computation
  sequence for $\varphi_e (x)$. \\(b) There is a {\tmem{universal predicate}}
  $\mathfrak{T_1}$ which satisfies $\forall e \forall x$:
  \[ P_e (x) \!\downarrow z \hspace{1em} \leftrightarrow \hspace{1em} \exists y
     \in 2^{\mathbbm{N}} [ \mathfrak{T_1} (e, x, y) \wedge Last (y) = z] . \]
\end{theorem}

Moreover as a corollary (to Theorem \ref{lzs-II}):

\begin{corollary}
  \label{jumps}(i) $x^{\nabla} \equiv_1 \Sigma_1$-$\tmop{Th} (\langle
  L_{\lambda^x} [x], \in, x \rangle)$  - the latter the $\Sigma_1$-theory of the
  structure.

  (ii) Let $x^{\infty} =_{\tmop{df}} \{e \mid \exists y P_e (x) \uparrow y\}$
  be the set of {\tmem{$x$-stable indices}}, of those program numbers whose output tapes eventually stabilize.
  Then $$\text{$x^{\infty} \equiv_1
  \Sigma_2$-$\tmop{Th} (\langle L_{\zeta^x} [x], \in, x \rangle)$}.$$
\end{corollary}

The conclusions are that the $\Sigma_{1}$-Theory of $L_{\lambda}$ is recursively isomorphic to the ITTM-jump $0^{\nabla}$.  One should compare this with Kleene's $\mathcal{O}$ being recursively isomorphic to the hyperjump, or again the $\Sigma_{1}$-Theory of $L_{\omega^{ck}_{1}}$.  The second part of the Corollary gives the analogous results for the index set of the eventually stable programs: here we characterise  $0^{\infty}$ as the $\Sigma_{2}$-Theory of $L_{\zeta}$.  The relativisations to inputs $x$ are immediate.

One should remark that extensions of Kleene's $\mathcal{O}$ from the ? case to the ITTM case are straightforward: we can define $\mathcal{O}^{+}$ by adding in to $\mathcal{O}$ those indices of Turing programs that now halt at some transfinite time.  After all, we are keeping the programs the same for both classes of machines, so we can keep the same formalism and definitions (literally) but just widen the class of what we consider computations. Similarly we can expand $\mathcal{O}^{+}$  to  $\mathcal{O}^{\infty}$ by adding in those indices of eventually stabilising programs.  This is done in detail in \cite{Kle07}.  We thus have:
$$ \frac{\mathcal{O}}{L_{\omega_{1}^{ck}}}  \approx
\frac{\mathcal{O}^{+}}{L_{\lambda}}  \approx
\frac{\mathcal{O}^{\infty}}{L_{\zeta}}.
$$

\subsection{Comparisons with Kleene recursion}

We have alluded to Kleene recursion in the introduction. His theory of recursion in higher types 
({\cite{Kl59}},{\cite{Kl62b}},{\cite{Kl62a}}, {\cite{Kl63}}) was an equational calculus, a
generalization of that for the G{\"o}del-Herbrand generalized recursive
functions. In this theory numbers were objects of type 0, whilst a function $f
: \mathbbm{N}^m \longrightarrow \mathbbm{N} $ is an object of type 1; and $F :
  \mathbbm{N}^l \times (2^{\mathbbm{N}})^m\longrightarrow \mathbbm{N}$ one of
type 2 etc. The $e${\tmem{'th procedure}} (whether thought of as the
$e$'th program of the kind of machine as outlined in the introduction, or else
as $e$'th equation system in his calculus) then allowed a computation with
inputs $\vec{n}, \vec{x}$ with oracle $\mathcal{I}$ to be presented in the
form $\{ e \} (\vec{n}, \vec{x}, \mathcal{I})$. The oracle $\mathcal{I}$ was
usually taken to include an oracle for existential quantification
$\mathcal{E}$ where, for $x \in 2^{\mathbbm{N}}$:

$$
\begin{array}{lcl}
& &  0 \mbox{ if } \exists n x (n) = 0 \\
\mathcal{E} (x) & = &  \\
& & 1 \mbox{ otherwise.}
\end{array}
$$

The reason for this was, although for any oracle $\mathcal{I}$ the class of
relations semi-decidable in $\mathcal{I}$ was closed under
$\forall^{\mathbbm{N}}$ quantification, when semi-decidable additionally in $\mathcal{E}$
it becomes closed under $\exists^{\mathbbm{N}}$ quantification. The Kleene
semi-decidable sets then would include the arithmetic sets in $\mathbbm{N}
\times 2^{\mathbbm{N}}$ (or further products thereof).  (Ensuring
computations be relative to $\mathcal{E}$ also guarantees that we have the
{\tmem{Ordinal Comparison Theorem}}.)

The decidable relations turn out to be the hyperarithmetic ones, and the
semi-decidable are those Kleene-reducible to $\tmop{WO}$, the latter being a
complete $\Pi^1_1$ set of reals. Thus:

\begin{theorem}
  \label{Kleene}{\em (Kleene)} The hyperarithmetic relations $R ( \text{$\vec{n}$,
  $\vec{x}$}) \subseteq \mathbbm{N}^k \times ( \mathbbm{N}^{\mathbbm{N}})^l$
  for any $k, l \in \mathbbm{N}$, are precisely those computable in
  $\mathcal{E}$.

  The $\Pi^1_1$ relations are precisely those semi-computable in
  $\mathcal{E}$.
\end{theorem}

Then a reducibility ordering comes from:

\begin{definition}
  \label{Kleenered}$\text{}${\tmstrong{Kleene reducibility:}} \ Let $A, B
  \subseteq \mathbbm{R}$; we say that
 $A$ is {\tmem{Kleene-semi-computable in $B$}} iff there is an index $e$ and
  $y \in \mathbbm{R}$ so that
$$\forall x \in \mathbbm{R}(x
  \in A \longleftrightarrow \{e\}(x, y, B, \mathcal{E}) \downarrow 1)) .$$
   $A$ is  {\tmem{Kleene computable in $B$}}, written, $A \leq_K B$,  iff 
  both $A$ and its complement are Kleene-semi-computable in $B$.
\end{definition}

The presence of the real $y$ deserves some explanation. We want to think of
the degree structure as being between sets of reals; the set $y$ throws in a
countable amount of information to the computation, and we are thus thinking
of two sets of reals $A =_K B$ as being of the same complexity up to this
countable amount of data. It implies that each Kleene degree contains
continuum many sets of reals, but moreover is closed under continuous
pre-images - it thus forms also a union of {\tmem{Wadge degrees}}.

We thus shall have that besides $\varnothing, \mathbbm{R}$ the bottommost
Kleene degree contains precisely all the Borel sets, whilst the degree of
$\tmop{WO}$ contains all co-analytic sets. As one sees the notion is very tied
up with hyperarithmeticity.

If we have a transitive reducibility notion $\leq$ on sets of integers $x$ say, together with a concomitant {\em jump operator} $x
\longrightarrow x'$ then an ordinal assignment $x \longrightarrow \tau^x\in On$ is
said to be a {\tmem{Spector criterion}} if we have:
$$x \leq y
\longrightarrow (x' \leq y \longleftrightarrow \tau^x < \tau^y). 
\ \ \ \ \ \ \ \ \ \ \ \ (*) $$

As an example if we take here {\tmem{hyperdegree}} $x \leq_h y$ (abbreviating
``$x$ is hyperarithmetic in $y$'') and the {\tmem{hyperjump}} operation, $x
\longrightarrow x^h$ where $x^h$ is a complete $\Pi^{1, x}_1$ set of integers,
then the assignment $x \longrightarrow \omega^x_{1 \tmop{ck}}$ \ (where the
latter is the least ordinal not (ordinary) Turing recursive in $x$) satisfies
the Spector Criterion $(*)$ above. \ \ For sets of reals $B$ we may extend this
notation and let $\omega^{B, x}_{1 \tmop{ck}}$ be the ordinal height
$\alpha$ of the least model of $\tmop{KP}$ set theory (so the least admissible
set) of the form $L_{\alpha} [x, B] \models \tmop{KP}$.

With this we may express $A \leq_KB$ as follows:

\begin{lemma}
  \label{Kdegree}$A \leq_K B$ iff there are $\Sigma_1$-formulae in
  $\mathcal{L}_{\in, \dot{X}} $ $\varphi_1 ( \dot{X}, v_0, v_1), \varphi_2 (
  \dot{X}, v_0, v_1)$, and there is $y \in \mathbbm{R}$, so that
\begin{center}
$\begin{array}{rcl}
\forall x \in
  \mathbbm{R} (x \in A & \Longleftrightarrow & L_{\omega_1^{B, y, x}} [B, y, x]
  \models \varphi_1 [B, y, x] \\ &  \Longleftrightarrow &
  L_{\omega_1^{B, y, x}} [B, y, x]
  \models \neg \varphi_2 [B, y, x]).
\end{array}
$
\end{center}
\end{lemma}

{\em Back to ITTM-semidecidability:}\\

The notion of semi-decidability comes in two forms.

\begin{definition}
  \label{reducibilities}(i) A set of integers $x$ is\tmtextit{{\tmem{
  semi-decidable}}} in a set y if and only if:
  \[ \exists e \forall n \in x \left[ P^y_e (n) \!\downarrow 1
     \longleftrightarrow n \in x \hspace{0.25em} \hspace{0.25em} \right] \]
  (ii) A set of integers $x$ is \tmtextit{{\tmem{decidable}} in a set }$y$ if
  and only if both $x$ and its complement is semi-decidable in $y$.
    We write $x \leq_{\infty} y$ for the reducibility ordering.\\
  (iii) A set of integers $x$ is
  \tmtextit{{\tmem{eventually-(semi)-decidable}}} in a set y if and only if
  the above holds with $\uparrow$ replacing $\!\downarrow$. For this
  reducibility ordering we write $x \leq_{}^{\infty} y$. \ \ \ \ \ \ \ \ \
  \ \ \
\end{definition}

We then get the analogue of the Spector criterion using $x^{\nabla}$ as the
jump operator:

\begin{lemma}
  \label{Spector}(i) The assignment $x \rightarrowtail \lambda^x$ satisfies
  the Spector Criterion:

  \ \ \ \ \ \ \ \ \ \ \ \ \ \ \ \ \ \ \ \ \ \ $x \leq_{\infty} y
  \longrightarrow (x^{\nabla} \leq_{\infty} y \leftrightarrow \lambda^x <
  \lambda^y)$.

  (ii) Similarly for the assignment $x \rightarrowtail \zeta^x :$

  \ \ \ \ \ \ \ \ \ \ \ \ \ \ \ \ \ \ \ \ \ \ $x \leq^{\infty} y
  \longrightarrow (x^{\infty} \leq^{\infty} y \leftrightarrow \zeta^x <
  \zeta^y)$
\end{lemma}

One can treat the above as confirmation that the ITTM degrees and jump
operation are more akin to hyperarithmetic degrees and the hyperjump, than to
the (standard) Turing degrees and Turing jump. Indeed they are intermediate
between hyperdergees and $\Delta^1_2$-degrees. 

To see this, we define a notion of degree using definability and Turing-invariant functions on reals (by the latter we mean a function
 $f : \mathbbm{R} \longrightarrow
\omega_1$ such that $x \leq_T y \longrightarrow f (x) \leq f
(y)$). 
Now assume that $f$ is
$\Sigma_1$-definable over $(\tmop{HC}, \in)$ without parameters, by a formula
in $\mathcal{L}_{\dot{\in}}$.

\begin{definition}
  \label{slices}Let $f$ be as described; let $\Phi$ be a class of formulae of
  $\mathcal{L}_{\dot{\in}}$. Then $\Gamma = \Gamma_{f, \Phi}$ is the
  pointclass of sets of reals $A$ so that $A \in \Gamma$ if and only if there
  is $\varphi \in \Phi$ with:

 $$\forall x \in
  \mathbbm{R} (x \in A \leftrightarrow L_{f (x)} [x] \models \varphi [x]).$$
\end{definition}

With  the function $f (x) = \omega^x_{1
\tmop{ck}}$ and $\Phi$ as the class of $\Sigma_{1}$-formulae we have that $\Gamma_{f,\Phi}$
coincides with the $\Pi^{1}_{1}$-sets of reals (by the Spector-Gandy Theorem). Replacing $f$ with the function $g(x)=\lambda^{x}$ then yields the (lightface) ITTM-semi-decidable sets. Lemma \ref{Kdegree} is then the relativisation of Kleene recursion which yields the relation $A\leq_{K} B$.

We now make the obvious definition:

\begin{definition}
  A set of reals $A$ is {\tmem{\tmtextit{{\tmem{semi-decidable}}} }}in a set
  of reals $B$ if and only if:
  \[ \exists e \forall x \in 2^{\mathbbm{N}} \left[ P^B_e (x) \downarrow 1
     \leftrightarrow x \in A \hspace{0.25em} \hspace{0.25em} \right] \]
  (ii) A set of reals $A$ is \tmtextit{{\tmem{decidable}} in a set of reals
  $B$} if and only if both $A$ and its complement is semi-decidable in $B$.\\
     (iii) If in the above we replace $\downarrow$ everywhere by $\uparrow$ then we obtain the notion in (i)  of {\em $A$ is eventually decidable in $B$} and in (ii) of {\em $A$ is eventually semi-decidable in $B$}.
\end{definition}

Then the following reducibility generalizes that of Kleene recursion.
\begin{definition}\mbox{ }\\
  (i) $A \leq_{\infty}\!B$  iff for some $e \in \omega$, for
  some $y \in 2^{\mathbbm{N}} : A$ is  decidable in $(y, B)$.\\
  (ii) $A \leq^{\infty}\!B$ iff for some $e \in \omega$, for
  some $y \in 2^{\mathbbm{N}} : A$ is eventually decidable in $(y, B)$.
\end{definition}

Again a real parameter has been included here in order to have degrees closed
under continuous pre-images. We should expect that  these reducibilities are
dependent on the ambient set theory, just as they are for Kleene degrees:
under $V = L$ there are many incomparable degrees below that of the complete
semi-decidable degree, and under sufficient determinacy there will be no
intermediate degrees between the latter and $\tmmathbf{0}$, and overall the
degrees will be wellordered. Now we get the promised analogy lifting Lemma \ref{Kdegree}, again generalizing in two ways depending on the reducibility.

\begin{lemma}
  \label{red-for-sets}\mbox{ } \\(i) $A \leq_{\infty}$B iff there are $\Sigma_1$-formulae
  in $\mathcal{L}_{\in, \dot{X}} $ $\varphi_1 ( \dot{X}, v_0, v_1), \varphi_2
  ( \dot{X}, v_0, v_1)$, and $y \in \mathbbm{R}$, so that

$  \begin{array}{rcl}
\forall x \in
  \mathbbm{R} (x \in A & \Longleftrightarrow & L_{\zeta^{B, y, x}} [B, y, x]
  \models \varphi_1 [B, y, x] \\
 & \Longleftrightarrow  & L_{\zeta^{B, y, x}} [B, y, x]
  \models \neg \varphi_2 [B, y, x]).
  \end{array}$

\noindent  (ii) $A \leq^{\infty}$B iff there are $\Sigma_2$-formulae in
  $\mathcal{L}_{\in, \dot{X}} $ $\varphi_1 ( \dot{X}, v_0, v_1), \varphi_2 (
  \dot{X}, v_0, v_1)$, and $y \in \mathbbm{R}$, so that

  $ \begin{array}{rcl}
  \forall x \in
  \mathbbm{R} (x \in A & \Longleftrightarrow & L_{\zeta^{B, y, x}} [B, y, x]
  \models \varphi_1 [B, y, x] \\
& \Longleftrightarrow & L_{\zeta^{B, y, x}} [B, y, x]
  \models \neg \varphi_2 [B, y, x]).
 \end{array}
  $
\end{lemma}

We have not formally defined all the terms here:  $\lambda^{B, y, x}$ is the supremum of the ordinals written by Turing programs acting transfinitely with oracles for $B,y$. The ordinal $\zeta^{B, y, x}$ is the least that is
not ITTM-$(B, x, y)$-eventually-semi-decidable.   There is a corresponding $\lambda$-$\zeta$-$\Sigma$-theorem and thus we have also that this $\zeta$ is 
 least such
that $L_{\zeta^{B, y, x}} [B, y, x]$ has a proper $\Sigma_2$-elementary
end-extension in the $L[B, y, x]$ hierarchy.

\subsection{Degree theory and complexity of ITTM computations}

Corollary \ref{jumps} shows that the ITTM-jump of a set of integers $x$ is
essentially a mastercode, or a $\Sigma_1$-truth set if you will, namely that of
$L_{\lambda^x} [x]$. \ The analogy here then is with $\mathcal{O}^x$, the
hyperjump of $x$, which is a complete $\Pi^{1,x}_1$ set of integers, \ as being
also recursively isomorphic to $\Sigma_1$-(Th$(L_{\omega^x_{1, \tmop{ck}}}
[x])$ . This again indicates that the degree analogy here should be pursued with
hyperdegrees rather than Turing degrees. It is possible to iterate the jump
hierachy through the $=_{\infty}$ - degrees, and one finds that, {\tmem{inside
$L$,}} the first $\zeta$-iterations form a linearly ordered hierachy with
least upper bounds at limit stages. We emphasise this as being inside $L$
since one can show that there is no least upper upper bound to $\{
0^{\nabla_n} \mid n < \omega \}$, but rather continuum many minimal upper
bounds. (see {\cite{W5}}). We don't itemize these results here but refer the reader instead to
{\cite{W}}.

A more general but basic open question is:

{\tmem{Q If $D =\{d_n : n < \omega\}$ is a countable set of
$=_{\infty}$-degrees, does $D$ have a minimal upper bound?}}

The background to this question is varied: for hyperdegrees this is also an
open question. Under Projective Determinacy a positive answer is known for
$\Delta^1_{2 n}$-degrees, but for $\Delta^1_{2 n + 1}$-degrees this is open,
even under PD. Minimal infinite time $\infty$-degrees can be shown to exist by similar
methods, using perfect set forcing, to those of Sacks for minimal hyperdegrees
(again see {\cite{W5}}).

One can also ask at this point about the nature of Post's problem for semi-decidable sets of integers. By the hyperdegree analogy one does not expect there to be incomparable such sets below $0^{\nabla}$ and indeed this turns out to be the case (\cite{HL2}).

\subsection{Truth and arithmetical quasi-inductive sets}

It is possible to relate ITTM's closely to an earlier notion due to Burgess {\cite{B}}
of {\tmem{arithmetical quasi-inductive definitions.}} We first make a general
definition:

\begin{definition}
Let  $\Phi : \mathcal{P}( \mathbbm{N}) \rightarrow \mathcal{P}( \mathbbm{N})$ be
  a $\Gamma$-operator, that is  ``$n\in \Phi(X)$'' is a $\Gamma$-relation. We define the $\Gamma$-{\tmem{quasi-inductive
  operator}} using iterates of $\Phi$ as:

  $\Phi_0 (X) = X ;$ \ \ \ \ \ $\Phi_{\alpha + 1} (X) = \Phi (\Phi_{\alpha}
  (X)) ;$

  $\Phi_{\lambda} (X) = \liminf_{\alpha \rightarrow \lambda} \Phi_{\alpha}
  (X)  =_{df} \cup_{\alpha < \lambda} \cap_{\lambda > \beta > \alpha} \Phi_{\beta}
  (X) .$

  We set the {\tmem{stability set}} to be $\Phi_{\tmop{On}} (X)$.
\end{definition}

By the nature of the $\liminf$ operation at limits, it is easy to see that the operation of an ITTM is an example of a recursive quasi-inductive operator on $\mathbbm{N}$.  Recall that a set of integers $B$ say is {\em inductive} if it is (1-1) reducible to the least fixed point of 
a monotone $\Pi^1_{1}$-operator. Burgess defined such a $B$ to be {\em arithmetically quasi-inductive} if it was (1-1) reducible to the {\em stability set} $\Phi_{\tmop{On}} (\varnothing)$.

In order to prove that an AQI halts, or reaches a stability point, one needs
to know that one has sufficiently long wellorderings, and a certain amount of
second order number theory is needed to prove that such ordinals exist. For
the case of the ITTM's we know which ordinals we need: $\Sigma^x$ for a
computation involving integers and the input real $x$. We then adopt this
idea of a `repeat pair' of ordinals for a quasi-inductive
operator $\Phi$: the least pair $(\zeta, \Sigma) =$ $(\zeta (\Phi, x),
\Sigma (\Phi, x))$ with $\Phi_{\zeta} (x) = \Phi_{\Sigma} (x) =
\Phi_{\tmop{On}} (x)$.

\begin{definition}
  $\mathsf{{AQI}}$ is the sentence: ``For every arithmetic operator
  $\Phi$, for every $x \subseteq \mathbbm{N},$ there is a wellordering $W$ with
  a repeat pair $(\zeta (\Phi, x), \Sigma (\Phi, x))$ in $\tmop{Field} (W)$''.
  If an arithmetic operator $\Phi$ acting on {\tmem{$x$}} has a repeat pair,
  we say that $\Phi$ {\em converges} (with input x).
\end{definition}

We may simulate an AQI with starting set $x\subseteq \mathbbm{N}$ as an ITTM with input $x$.
Since we know how long such  a machine takes to halt or loop, this gives the length of ordering needed to determine the extent of the AQI. Given the characterisation from the (relativized) $\lambda$-$\zeta$-$\Sigma$-Theorem one arrives at the fact that
\begin{theorem} The theories
  $\mathsf{\Pi^{1}_3}$-$\mathsf{{CA}_0}$,
  $\mathsf{\Delta^{1}_3}$-$\mathsf{{CA}_0}{+ \mathsf{AQI}}$,
  and $\mathsf{\Delta^{1}_3}$-$\mathsf{{CA}_0}$ are in
  strictly descending order of strength, meaning that each theory proves the
  existence of a $\beta$-model of the next.
\end{theorem}

What was engaging Burgess was an analysis of a {\tmem{theory of truth}} due to
Herzberger {\cite{H}}. \ The latter had defined a {\tmem{Revision
Sequence}} which was essentially a quasi-inductive operator, just a bit beyond
the arithmetic as follows.

$H_{0} =\emptyset$:

 $H_{\alpha + 1} =\{ \ulcorner \sigma \urcorner : \langle \mathbbm{N}, +,
\times, \cdots, H_\alpha \rangle \models \sigma\}$ ; with $H_\alpha$ interpreting $T$;

 $H_{\lambda} = \cup_{\alpha < \lambda} \cap_{\lambda > \beta > \alpha}
H_{\beta}$.

Burgess then defined the AQI sets as above and calculated that the ordinals
$(\zeta, \Sigma)$ formed exactly the repeat pair needed for AQI's or for the Herzberger
revision sequence. This was much earlier than the invention of ITTM's and was
unknown to workers in the latter area around 2000, until Benedikt L\"owe pointed
out ({\cite{L}}) the similarity between the Herzberger revision sequence formalism and that
of the machines. 
It can be easily seen that any
Herzberger sequence with starting distribution of truth values $x$ say, can be
mimicked on an ITTM with input $x$.  Thus this is one way of seeing that
Herzberger sequences must have a stability pair lexicographically no later
than $(\zeta, \Sigma)$. Burgess had shown that H-sequences then loop at no
earlier pair of points.  More recently Field {\cite{Fi03}} has used a revision
theoretic definition with a $\Pi^1_1$-quasi-inductive operator to define a
variant theory of truth.  For all three formalisms, Fields, Burgess's AQI, and
ITTM's, although differing considerably in theory, the  operators are all
essentially equivalent as is shown in {\cite{W16}}, since they produce
recursively isomorphic stable sets. The moral to be drawn from this is that in
essence the strength of the liminf rule is at play here, and seems to swamp all else.

\section{Variant ITTM models}

Several questions readily occur once one has formulated the ITTM model. Were
any features chosen crucial to the resulting class of computable functions? \
Do variant machines produce different classes? Is it necessary to have three
tapes in the machine? \ The answer for the latter question is both yes and no. First the
affirmative part: it was shown in {\cite{HaSe}} the class of functions $f :
\mathbbm{N} \longrightarrow \mathbbm{N}$ remains the same if 3 tapes are
replaced by 1, but not the class of functions $f : 2^{\mathbbm{N}}
\longrightarrow 2^{\mathbbm{N}}$. \ The difficulty is somewhat arcane: one may
simulate a 3-tape machine on a 1-tape machine, but to finally produce the
output on the single tape and halt, some device is needed to tell the machine
when to finish compacting the result down on the single tape, and they show
that this cannot be coded on a 1-tape machine. On the other hand {\cite{W8}}
shows that if one adopts an alphabet of three symbols this can be done and the
class of functions $f : 2^{\mathbbm{N}} \longrightarrow 2^{\mathbbm{N}}$ is
then the same. One may also consider a B for `Blank' as the third  symbol, and change
the liminf rule so that if cell $C_i$ has varied cofinally in a limit ordinal
$\lambda$, then $C_i (\lambda)$ is set to be blank (thus nodding towards
ambiguity of the cell value).  With this alphabet and liminf
rule a 1-tape machine computes the same classes as a 3-tape machine, and these are both the same as computed by the original ITTM.

What of the liminf rule itself? We have just mentioned a variant in the last
paragraph. Our original liminf rule is essentially of a $\Sigma_2$ nature: a
value of 1 is in a cell $C_i (\mu)$ at limit time $\mu$ if there is an $\alpha
< \mu$ such that for all $\beta \in (\alpha, \beta)$ $C_i (\beta) = 1$. \ \
Running a machine inside $L_{\mu}$ one sees that the snapshot $s (\mu)$ is a
predicate that is $\Sigma_2$-definable over $L_{\mu}$.  It was observed in {\cite{W}}  that
the liminf rule is {\tmem{complete}} for all other rules $\Sigma_2$-definable over
limit levels $L_{\mu}$ in that for any other such rule the
stability set obtained for the universal machine (on 0 input) with such a rule is (1-1)
  $\Sigma_2$-definable over
$L_{\zeta}$ and thus is (1-1) in the $\Sigma_2$-truth set for $L_{\zeta}$. 
However the latter is recursively isomorphic to the stability set for the universal ITTM by Corollary
\ref{jumps} and hence the standard stability set subsumes that of another machine with a different $\Sigma_{2}$-rule.  Given the $\Sigma_2$ nature of the limit rule, with hindsight one
sees that it is obvious that with $(\zeta', \Sigma')$  defined to be the
lexicographically least pair with $L_{\zeta'} \prec_{\Sigma_2} L_{\Sigma'}$,
then we must have that the universal ITTM  enters a loop at $\zeta'$. That
it cannot enter earlier of course is the $\lambda$-$\zeta$-$\Sigma$-Theorem, but a vivid way to see
that this is the case is afforded by the construction in {\cite{FrWe0?}} which
demonstrated that there was a non-halting ITTM program producing on its
output tape continually sets of integers that coded levels $L_{\alpha}$ of the
constructible hierarchy for ever larger $\alpha$ below $\Sigma$; at stage
$\Sigma$ it would perforce produce the code for $L_{\zeta}$ and then forever
cycle round this loop producing codes for levels $\alpha \in [\zeta, \Sigma)$.

More complex rules lead to more complex machines. These were dubbed
`hypermachines' in {\cite{FrWe11}}, where a machine was defined with a
$\Sigma_3$-limit rule, and this was shown to be able to compute codes for
$L_{\alpha}$ for $\alpha < \Sigma (3)$, where now $\zeta (3) < \Sigma (3)$
was the lexicographically least pair with $L_{\zeta (3)} \prec_{\Sigma_3}
L_{\Sigma (3)}$. The stability set was now that from the snapshot at stage
$\zeta (3)$ and was (1-1) to the $\Sigma_3$-truth set for this level of $L$.
Inductively then one defines $\Sigma_4, \Sigma_5, \ldots \nocomma, \Sigma_n,
\ldots$ limit rules with the analogous properties. I think it has to be said
though that the definitions become increasingly complex and even for $n = 3$,
mirror more the structure of $L$ in these regions with its own `stable
ordinals' rather than anything machine-inspired. With these constructions one
can then `compute' any real that is in $L_{\tau}$ where $\tau = \sup_n \zeta
(n)$.

\subsection{Longer tapes}

generalizations of the ITTM machine are possible in different directions. One
can consider machines with tapes not of cells of order type $\omega$ but of
longer types. Some modifications are needed: what do we do if the program asks
the R/W head to move one step leftwards when hovering over a cell
$C_{\lambda}$ for $\lambda$ a limit ordinal? \ There are some inessentially
different choices to be made which we do not catalogue here, but assume some
fixed choices have been made.

We consider first the extreme possibility that the tape is of length
$\tmop{On}$, that is of the class of all ordinals. We now have the possibility
that arbitrary sets may be computed by such machines.  Independently Dawson and Koepke came up with this concept. There are some caveats:
how do we know that we can `code' sets by transfinite strings of $0, 1$'s at
all? Dawson
({\cite{Da09}}) formulated an {\tmem{Axiom of Computability}} that said every
set could appear coded on the output tape of such a machine at some stage
whilst it was running; thus for any set $z$ there would be a program number
$e$  with $P_e$ (not necessarily halting) with a code for $z$ appearing on
its output tape. He then argued that the class of such sets is a model of ZFC,
and by studying the two dimensional grid of snapshots produced a
L\"owenheim-Skolem type argument to justify that the Axiom of Computability
implied the Generalized Continuum Hypothesis. \ That the class of computable
sets satisfied $\tmop{AC}$ falls out of the assumption that sets can be coded
by strings and such can be ordered. Since this machine's operations are again
very absolute, it may be run inside $L$, thus demonstrating that `computable
sets' are nothing other than the constructible sets. \ Koepke in {\cite{K05}}
and later with Koerwien in {\cite{KK06}} considered instead halting
computations starting with an $\tmop{On}$-length tape marked with finitely
many $1$'s in certain ordinal positions $(n,\xi
_{1}, \ldots, \xi_{n} )$,  and asked for a computation as to
whether $(\varphi_{n}(\xi_1, \ldots \nocomma \xi_{n - 1}))^{L_{\xi_n}}$ was true. Thus the machine was capable of computing a truth predicate for
$L$. This leads to:

\begin{theorem}
  {\tmem{(Koepke {\cite{K05}})}} A set $x \subseteq \tmop{On}$ is
  $\tmop{On}$-$\tmop{ITTM}$-computable from a finite set of ordinal parameters
  if and only if it is a member of the constructible hierarchy.
\end{theorem}

One might well ask whether the computational approach to $L$ might lead to
some new proofs, or at least new information, on some of the deeper fine
structural and combinatorial properties of $L$. However this hope turned out
to be seemingly thwarted by the $\Sigma_2$-nature of the limit rule. Fine
structural arguments are very sensitive to definability issues, and in
constructions such as that for Jensen's $\Box$ principle, say, we need to know when or how ordinals are
singularised for any $n$ including $n = 1$ and the limit rule works against
this. Moreover alternatives such as the {\tmem{Silver Machine}} model which
was specifically designed to by-pass Jensen's fine structural analysis of $L$,
make heavy use of a {\tmem{Finiteness Property}} that everything appearing at
a successor stage can be defined from the previous stages and a finite set of
parameters;  just does not seem to work for $\tmop{On}$-ITTM's.

However this does bring to the fore the question of shortening the tapes to
some admissible ordinal length $\alpha > \omega$ say, and asking what are the
relations between $\alpha$-ITTM's and the $\alpha$-recursion theory developed
in the late 1960's and early 70's. The definitions of that theory included that a
set $A \subseteq \alpha$ which is $\Sigma_1 (L_{\alpha})$ was called
$\alpha${\tmem{-recursively enumerable}} ($\alpha$-r.e.). It was
$\alpha$-{\tmem{recursive}} if both it and its complement is $\alpha$-r.e. and
thus is $\Delta_1 (L_{\alpha})$.  A notion of  {\tmem{relative $\alpha$-recursion}} was defined but then noticed to be intransitive; a stronger notion was defined and  denoted by $A
\leq_{\alpha} B$.

Koepke and Seyfferth in {\cite{KoSe09}} define $A${\tmem{ is computable in
$B$}} to mean that the characteristic function of $A$ can be computed by a
machine in $\alpha$ many stages from an oracle for $B$. This is exactly the
relation that $A \in \Delta_1 (L_{\alpha} [B])$. This has the
advantage that the notion of $\alpha$-computability and the associated
$\alpha$-computable enumerability ($\alpha$-c.e.) tie up exactly with the notions of
$\alpha$-recursiveness and $\alpha$-r.e. They then reprove the Sacks-Simpson
theorem for solving Post's problem: namely that there can be two $\alpha$-c.e.
sets neither of which are mutually computable in their sense from the other.

However the relation ``{\tmem{is computable in}}'' again suffers from being an intransitive
one. Dawson defines the notion of $\alpha$-{\tmem{sequential computation}}
that requires the output to the $\alpha$-length tape be written in sequence
without revisions. This gives him a transitive notion of relative
computability: a set is $\alpha$-computable if and only if it is
$\alpha$-recursive, and it is $\alpha$-computably enumerable if and only if
it is both $\alpha$-r.e. and {\em regular}. Since Sacks had shown ({\cite{Sa66}})
that any $\alpha$-degree of $\alpha$.-r.e. sets contains a regular set, he
then has that the structure of the $\alpha$-degrees of the $\alpha$-r.e. sets
in the classical, former, sense, is isomorphic to that of the $\alpha$-degrees of the
$\alpha$-c.e. sets. This implies that theorems of classical $\alpha$-recursion
theory about $\alpha$-r.e. sets whose proofs rely on, or use regular
$\alpha$-r.e. sets will carry over to his theory. This includes the
Sacks-Simpson result alluded to. The Shore Splitting Theorem ({\cite{Sho75}})
which states that any regular $\alpha$-r.e. set $A$ may be split into two
disjoint $\alpha$-r.e. sets $B_0, B_1$ with $A \nleq_{\alpha} B_i$,  is less
amenable to this kind of argument but with some work the Shore Density theorem
({\cite{Sho76}}) that between any two $\alpha$-r.e. sets $A <_{\alpha}
B$ there lies a third $\alpha$-r.e. $C$: $A <_{\alpha} C <_{\alpha} B$ can be
achieved. As Sacks states in his book, the latter proof seems more bound up with the
finer structure of the constructible sets than the other $\alpha$-recursion
theory proofs.  Dawson generalizes this by lifting his notion of
$\alpha$-computation to that of a $\mathbbm{B}$-$\alpha$-computation where now
$\underline{\mathbbm{B}}  =_{\tmop{df}} \langle J_{\alpha}^{\mathbbm{B}}, \in,
\mathbbm{B} \rangle$ is an admissible, {\tmem{acceptable}}, and {\tmem{sound}}
structure for a $\mathbbm{B} \subseteq \alpha$. These assumptions make
$J_{\alpha}^{\mathbbm{B}}$ sufficiently $L$-like to rework the Shore argument
to obtain:

\begin{theorem}
  {\tmem{(Dawson - The $\alpha$-c.e.~Density Theorem)}} Let $\text{$\underline{\mathbbm{B}}$}$ be as above. Let
  $A, B$ be two $\mathbbm{B}$-$\alpha$-c.e. sets, with $A <_{\mathbbm{B},
  \alpha} B$. Then there is $C$ also $\mathbbm{B}$-$\alpha$-c.e, with $A
  <_{\mathbbm{B}, \alpha} C <_{\mathbbm{B}, \alpha} B$.
\end{theorem}

\section{Other Transfinite Machines}

Once the step has been taken to investigate ITTM's, one starts looking at other machine models and sending them into the transfinite. We look here at {\tmem{Infinite Time Register Machines}}
(ITRM's) both with integer and ordinal registers, and lastly comment on 
{\tmem{Infinite Time Blum-Shub-Smale Machines}} (IBSSM's).

\subsection{Infinite time register machines (ITRM's)}

A (standard) {\tmem{register machine}} as devised by  Shepherdson and Sturgis
{\cite{SheStu63}}, or Minsky \ {\cite{Min67}}, consists of finite number of
natural number registers $R_i$ for $i < N$, running under a program consisting
of a finite list of instructions $\vec{I} = I_0, \ldots, I_q$. The latter
consist of zeroising, transferring of register contents one to another, or
conditional jump  to an instruction number in the
program, when comparing two registers.  At time $\alpha$ we shall list the $N$-vector of register contents as
$\vec{R} (\alpha)$. The next instruction the machine is about to perform we
shall denote by $I (\alpha)$.  We adopt a liminf rule again. Thus the next
instruction to be performed at limit stage $\lambda$, is $I (\lambda)
=_{\tmop{df}} \liminf_{\alpha \rightarrow \lambda} I (\alpha)$. As discussed before for ITTM's,
this can be seen to place control at the beginning of the outermost loop, or
subroutine, entered cofinally often before stage $\lambda$. We shall use a
$\liminf^{\ast}$ rule for register contents: if a register's contents edges up
to infinity at time $\lambda$ it is reset to $0$:
$$R_i (\lambda) =_{\tmop{df}} \liminf^{}_{\alpha \rightarrow \lambda} R_i
(\alpha) {\mbox{ if this is finite; otherwise we set }} R_i (\lambda) = 0.$$

Athough perhaps not apparent at this point, it is this `resetting to zero' that gives
the ITRM its surprising strength: specifying that the machine, or program,
crash with no output if a register becomes unbounded results in a
substantially smaller class of computable functions.  A function $F :
\mathbbm{N}^N \rightarrow \mathbbm{N}$ is then {\tmem{ITRM-computable}} if
there is an ITRM program $P$ with $P (\vec{k}) \da F (\vec{k})$ for
every $\vec{k} \in \mathbbm{N}^N$. In order to accommodate computation from
a set of integers $Z \subseteq \mathbbm{N}$ say, we add an oracle query
instruction $? k \in Z$? and receive as $0 / 1$ the answer to a register as a
result.

These machines were defined by Koepke and investigated by him and co-workers
in {\cite{CFKMNW}}, {\cite{KoMi08}}. A {\tmem{clockable ordinal}} has the same
meaning here as for ITTM's, except that here these ordinals for man initial
segment of $\tmop{On}$. Defining a {\tmem{computable ordinal}} as one which
has real code whose characteristic function is ITRM-computable, they show that
the clockables ordinals coincide with the computable ordinals. To analyse
what they are capable of, first note as a crude upper bound that they could be
easily simulated on an ITTM.
However ITRM's can detect whether an oracle set $Z \subseteq
\mathbbm{N}$ codes a wellfounded relation: a backtracking algorithm that
searches for leftmost paths can be programmed. Thus $\Pi^1_1$-sets are ITRM-decidable.

It also turns out that, in contradistinction to the finite case, the strength
of the infinite version of register machines  diverges from that of the Turing machine, but moreover there is no universal
ITRM.  We outline the arguments for this.
\begin{definition}
 $(N$-{\tmem{register halting set}}$)$
$$H_N
  =_{\tmop{df}} \{ \langle e, r_0, \ldots, r_{N - 1} \rangle \mid P^{}_e (r_0,
  \ldots, r_{N - 1}) \downarrow\}.$$
\end{definition}
(There is an obvious generalization $H_n^Z$ for machines with oracle $Z$.)

Koepke and Miller show that if there is some instruction $I'$ and register
contents vector $\overrightarrow{R_i}$ such that the snapshot $(I',
\overrightarrow{R_i})$ reoccurs in the course of computation at least
$\omega^{\omega}$ times in order type, then the computation is in a loop and
will go on for ever. 

\begin{theorem}
  \label{h-sets-exist} {\tmem{{\tmem{{\tmem{(Koepke-Miller {\cite{KoMi08}}) \
  }}For any $N$ the $N$-halting problem}}}}: `$\langle e, \vec{r} \rangle \in
  H_N$' is decidable by an ITRM. Similarly for any oracle $Z$, the $(N,
  Z)$-{\tmem{halting problem}} `$\langle e, \vec{r} \rangle \in H^Z_N$' is
  decidable by a $Z$-ITRM with an oracle for $Z$.
\end{theorem}

The number of registers has to be increased to calculate $H_N$ for
increasing $N$. The corollary to this is that there can be no one single
universal ITRM. We can get an exact description of the strength of ITRM's by
assessing bounds on the ordinals needed to see that a machine either halts or
is looping. \ It is discussed in {\cite{KoWe11}} and shown there that if one
has an ITRM with a single register then it has either halted or is in an
infinite loop by the second admissible ordinal $\omega^{\tmop{ck}}_2$. \ One
cannot replace this with $\omega_1 = \omega^{\tmop{ck}}_1$: if
$\tmop{Liminf}_{\beta \rightarrow \omega_1} R_0 (\beta) = p < \omega$ then a
$\Pi_2$-reflection argument shows that the same instruction number is used, and the
value in the register is this $p$, on a set of ordinals closed and unbounded
in $\omega_1$. By the Koepke-Miller criterion mentioned above this would
indeed mean that the computation was looping. However it can be the case that
$\tmop{Liminf}_{\beta \rightarrow \omega_1} R_0 (\beta) = \omega$ and then
this would have to be reset to 0: $R_0 (\omega_1) = 0$. Then again
$\tmop{Liminf}_{\beta \rightarrow \omega_1 + \omega_1} R_0 (\beta)$ may also
be $\omega$, but the instruction number could now differ. However by
$\omega^{\tmop{ck}}_2$ the criterion will have already applied and the
computation if still running will be looping. One then shows by induction that
each extra register added to the architecture requires a further admissible
ordinal in run time to guarantee looping behaviour. \ One then thus arrives at
the property that any ordinal below $\omega^{\tmop{ck}}_{\omega}$- the first
limit of admissibles, is clockable by such an ITRM, and thence that the
halting sets $H_n$ can be computed by a large enough device. We can state this more formally:

Thus the assertion that that these machines either halt or exhibit looping
behaviour turns out to be equivalent to a well known subsystem of second order
number theory, namely,   $\mathsf{\Pi^1_1}$-$\mathsf{CA_0}$. 
Let  $\mathsf{ITRM}_N$ be the assertion:
  {\noindent}``The $N$-\tmem{register halting set} $H_N$ exists.'' Further, let  {\tmstrong{{$\mathsf{ITRM}$}}} be
the similar relativized statement that
  ``For any $Z \subseteq \omega$, for any $N < \omega$ the
  $N$-{\tmem{register halting set}} $H^Z_N$ exists.''  Then more
precisely:

\begin{theorem}{\em (Koepke-Welch \cite{KoWe11})}\mbox{ }\\
   (i)  $\mathsf{\Pi^1_1}$-$\mathsf{CA_0}\vdash
{\tmstrong{{$\mathsf{ITRM}$}}}$.  In particular:\\
  $\mathsf{KP} +$``there exist $N+1$
  admissible ordinals $> \omega$'' $ \vdash \mathsf{ITRM}_N$.

\noindent  (ii) $\mathsf{ATR_0} + {{\tmstrong{{$\mathsf{ITRM}$}}}} \vdash
  \mathsf{\Pi^1_1}$-$\mathsf{CA_0}$.

   In particular there is a fixed $k < \omega$ so that for any $N < \omega$
 $$\mathsf{ATR_0} + \mathsf{ITRM}_{N \cdot k}
  \vdash\mbox{``}\tmop{HJ} (N, \varnothing) \mbox{ exists.'' }$$
\end{theorem}

 An analysis of Post's problem in this context is effected in \cite{HaMi09}.

{\subsection{ Ordinal register machines ORM's}}

We mention finally here the notion studied by Koepke and Siders of
{\tmem{Ordinal Register Machines}} (ORM's, {\cite{KoSi06}}): essentially these are the
devices above but extended to have ordinal valued registers.  Platek (in
private correspondence) indicated that he had originally considered his equational
calculus on recursive ordinals as being implementable on some kind
of ordinal register machine. Siders also had been thinking of such machines and
in a series of papers with Koepke considered the unbounded ordinal model. The
resetting $\tmop{Liminf}^{\ast}$ rule is abandoned, and natural
$\tmop{Liminf}$'s are taken. Now ordinal arithmetic can be performed.
Remarkably given the paucity of resources apparently available one has the
similar theorem to that of the On-ITTM:\\

\begin{theorem}
  {\tmem{\label{ORM-L}(Koepke-Siders {\cite{KoSi06}})}} A set $x \subseteq
  \tmop{On}$ is $\tmop{ORM}$-computable from a finite set of ordinals
  parameters if and only if it is a member of the constructible hierarchy.
\end{theorem}

They implement an algorithm that computes the truth predicate $T \subseteq
\tmop{On}$ for $L$ and which is ORM-computable on a 12 register machine (even
remarking that this can be reduced to 4!). From $T$ a class of sets
$\mathcal{S}$ can be computed which is a model of their theory $\tmop{SO}$,
which is indeed the constructible hierarchy.

\section{Infinite Time Blum-Shub-Smale machines IBSSM}

Lastly we consider the possible transfinite versions of the Blum-Shub-Smale
machines. These can be viewed as having registers $R_1 \nocomma, \ldots
\nocomma, R_N$  containing now euclidean reals $r_1, \ldots, r_n \in \mathbbm{R}$. \
There is a finite program or flow-chart with instructions divided into function nodes or conditional
branching nodes. We shall assume that function nodes have the possibility of
applying a rational function computation of the registers (we test each time
that we are not dividing by zero). So far this accords with the finite BSS
version.  We now make the, rather stringent, condition that at a limit stage
$\lambda$ if any register $R_i$ does not converge to a limit in the usual
sense, then the whole computation is deemed to have crashed and so be
undefined. The value of $R_{\iota} (\lambda)$ is then set to be the
ordinary limit of the contents of $R_i (\alpha)$ as $\alpha \rightarrow
\lambda$. It has been noted that a BSS machine cannot calculate the functions
$e^x$, $\sin x$ etc., but an IBSSM can, indeed in $\omega$ many steps (by simply calculating increasingly long initial segments of the appropriate power series). 

Koepke and Seyfferth {\cite{KoSe09}} have investigated such machines with
continuous limits. To simulate other sorts of machines on an IBSSM requires
some ingenuity: a register that is perhaps simulating a register of one of the
ITRM's discussed earlier, may have some contents $x$, that tends to infinity
and be then reset. Here then it is better to calculate with $\frac{1}{x}$ in
order to ensure a continuous limit of $0$. \ Else if the register is
simulating the contents of the scratch tape of an ITTM, then perhaps at
successor stages continual division by 2 ensures again a continuous limit at
the next limit ordinal of time. They show that a machine with $n$ nodes in its
flow diagram can halt on rational number input at ordinal times without any
limit below $\omega^{n + 1}$. Thus any such machine will halt, crash or be
looping by time $\omega^{\omega}$.

The question is naturally what is the computational power of such machines?
Clearly, by absoluteness considerations, on rational input such a machine can
be run inside the constructible hierarchy, and indeed from what they showed on
ordinal lengths of computations, inside $L_{\omega^{\omega}}$.  They then naturally ask
whether any real in $L_{\omega^{\omega}}$ can be produced by an IBSSM machine?

We can answer this affirmatively below. However at the same time we combine this with yet
another characterisation.  It is possible to give another characterisation of
the sets in $L_{\omega^{\omega}}$ by using the notions of {\tmem{Safe
Recursive Set Functions}} (SRSF) of Beckmann, Buss and S. Friedman
{\cite{BeBuFr12}}. They are generalizing  the notion of safe recursion of
Bellantoni and Cook (\cite{BellCook92}) used  to define polynomial time computations. \ Here
variables are divided into two types {\tmem{safe}} and{\tmem{ normal}}. In the
notation $f (\vec{a} / \vec{b})$ recursion is only allowed on the
safe variables in $\vec{b}$. This allows for the definition by recursion of
addition and multiplication but crucially not exponentiation. One of the aims
of {\cite{BeBuFr12}}  is to have a notion of set recursion that corresponds
to `polynomial time'.  On input an $\omega$-string in $2^{\mathbbm{N}}$
for example, one wants a computation that halts in polynomial time from
$\omega$ - the length of the input. Hence the calculation should halt by some
$\omega^n$ for an $n < \omega$. They have:

\begin{theorem} {\em({\cite{BeBuFr12}})}
Let $f$ be any $SRSF$. Then there is a ordinal polynomial $q_f$ in variables $\vec \alpha$ so that
$$ rk(f(\vec{a}/ \vec{b})) \leq \tmop{max}_i rk(a_i)  + q_f(\vec{rk({a}})).$$
\end{theorem}

Thus the typing of the variables ensures that the ranks of sets computed as
outputs from an application of an SRSFunction are polynomially bounded in the
ranks of the input. Using an adaptation of Arai, such functions on finite
strings correspond to polynomial time functions in the ordinary sense. For
$\omega$-strings we have that such computations halt by a time polynomial in
$\omega$. As mentioned by Schindler, it is natural to define `polynomial time'
for ITTM's to be those calculations that halt by stage $\omega^{\omega}$, and
a polynomial time ITTM function to be one that,  for some $N < \omega$, terminates on all inputs by
time $\omega^N$. We thus have:

\begin{theorem}
  The following classes of functions of the form $F : (2^{\mathbbm{N}})^k
  \rightarrow 2^{\mathbbm{N}}$ are extensionally equivalent:

  {\em (I)} Those functions computed by a continuous {\tmem{IBSSM}} machine;
  
  {\em (II)} Those functions that are polynomial time {\tmem{ITTM}};
  
   {\em (III)} Those functions that are safe recursive set functions.
\end{theorem}

\noindent Proof: We take $k = 1$. We just sketch the
ideas and the reader may fill in the details. \ By Koepke-Seyfferth for any
IBSSM computable function there is $N < \omega$ so that the function is
computable in less than $\omega^{N }$ steps. We may thus consider that
computation to be performed inside $L_{\omega^{N }} [x]$ and so potentially
simulable in polynomial time (in $\omega^{M}$ steps, for some $M$) by an ITTM. However
this can be realised:  a code for any $L_{\alpha} [x]$ for $\alpha \leq
\omega^N$, $x \in 2^{\mathbbm{N}}$, and its theory, may be computable by an
ITTM (uniformly in the input $x$) by time $\omega^{N + 3}$ by the argument of
Lemma 2 of {\cite{FrWe0?}} (Friedman-Welch).  Since we have the theory, we have the digits of the final halting IBSSM-output (or otherwise the fact that it is looping or has crashed respectively, since these are also part of the set theoretical truths of $L_{\omega^{N }} [x]$). Thus (II) $\supseteq$ (I). 
If $F$ is in the class (II), then for some $N < \omega$, $F (x)$ is computable
within $L_{\omega^N} [x]$ and by setting up the definition of the ITTM program
$P$ computing $F$ we may define some $\alpha$ such that the output of that
program $P$ on $x$ ({\em i.e.} $F (x)$) is the $\alpha$'th element of
$2^{\mathbbm{N}}$ in $L_{\omega^N} [x]$ \ uniformly in $x$. However the set
$L_{\omega^N} [x]$ is SRSF-recursive from $\omega \cup \{x\}$ (again uniformly
in $x$) as is a code for $\alpha$. This yields the conclusion that we may find
uniformly the output of $P (x)$  using the code for
$\alpha$, again as the output of an SRSF-recursive-in-$x$ function. \ This
renders (II) $\subseteq$ (III).

Finally if $F$ is in (III), (and we shall assume that the variable $x$ is in a
safe variable place - but actually the case where there are normal and safe variables is
handled no differently here) then there is (\tmtextit{cf.
}{\cite{BeBuFr12}}, 3.5) a finite $N$ and a $\Sigma_1$-formula
$\varphi (v_0, v_1)$ so that $F (x) = z$ iff $L_{\omega^N } [x] \models
\varphi [x, z]$ \ (using here that $\tmop{TC} (x) = \omega$ and thus
$\tmop{rk} (x) = \omega$). Indeed we may assume that $z$ is named by the
canonical $\Sigma_1$-skolem function $h$ for, say, $L_{\omega^N + \omega} [x]$ as $h
(i, n)$ for some $n < \omega$. Putting this together we have some $\Sigma_1
\,\psi (v_0)$ (in $\mathcal{L}_{\dot{x,} \dot{\in}}$) so that $F (x) (k) = z (k)
= 1$ iff $L_{\omega^N + \omega} [x] \models \psi [k]$. In short to be able to
determine such $F (x)$ by an IBSSM it suffices to be able to compute the
$\Sigma_1$-truth sets for $L_{\alpha} [x]$ for all
$\alpha < \omega^{\omega}$ by IBSSM's. There are a variety of ways one could
do this, but it is well known that calculating the $\alpha$'th iterates of the Turing jump relativised to $x$ for $\alpha<\omega^{\omega}$ 
would suffice. To simplify notation we shall let $x$ also denote the set of integers in the infinite fractional expansion of the real $x$.  So fix a $k<\omega$, to see that we may calculate $x^{(\beta)'}$ for $\beta < \omega^k$.
One first constructs a counter  to be used in general iterative processes, 
using registers $C_0, \ldots \nocomma
C_{k - 1}$  say, whose contents represent the integer coefficients in the Cantor normal form of $\beta<\omega^{k}$ where we are at the $\beta$'th stage in the process. (The counter of course must conform to the requirement that
registers are continuous at limits $\lambda \leq \omega^k$. This can be
devised using reciprocals and repeated division by $2$ rather than
incrementation by $1$ each time.) We assume this has
been done so that in particular that $C_0 = C_1 = \cdots = C_{k -
1} = 0$ occurs first at stage $\omega^k$. We then
  code the
characteristic function of $
\{ m \in \omega \mid m \in W_m^{x^{(\beta)'}} \}$ as $1 /
0$'s in the digits at the $s$'th-places after the decimal point of $R_1$ where
$s$ is of the form $p_{k+m}.p^{n_0+1}_{ 0} . \cdots .p^{n_{k - 1}+1}_{k - 1}$ \ where
$p_0 = 2,  p_1 = 3$, etc., enumerates the primes, and $n_{j}$ the exponent of $\omega^{j}$ in the Cantor Normal form of $\omega^{\beta}$.  For limit stages $\lambda<\omega^k$, continuity of the register contents automatically ensures that this real in $R_{1}$ also codes the disjoint union of the $
x^{(\beta)'}$ for $\beta < \lambda$, and at stage $\omega^{k}$ we have the whole sequence of jumps encoded as required. \mbox{ }\hfill \qed

\section{Conclusions}

The avenues of generalization of the Turing machine model into the
transfinite which we have surveyed, give rise to differing perspectives and a
wealth of connections. Higher type recursion theory, to which the models
mostly nearly approximate, to a lesser or greater extent, was a product of
Kleene's generalization of the notion of an equational calculus approach to
recursive functions.  Here discussed are machines 
more on the Turing side of the balance. Some of the other generalizations
of recursion theory, say to meta-recursion theory, as advocated by Kreisel and
elucidated by Sacks and his school, and which later became ordinal
$\alpha$-recursion theory, we have not really discussed here in great
detail, but again their motivations came from the recursion theoretic-side,
rather than any `computational-model-theoretic' direction. The
models discussed in this chapter thus fill a gap in our thinking.

Referring to the last section, we find that, rather like a Church's
thesis, we have here an effective system for handling  $\omega$-strings in
polynomial time, as formalized by the SRSF's, and a natural corresponding computational
model of ITTM's working with calculations halting by time earlier than
$\omega^{\omega}$. The model of computation with the continuous limit
IBSSM's then also computes the same functions. Note that assertions such as
 that ``{\tmem{every continuous IBSSM halts, loops, or becomes discontinuous}}'' when
formalized in second order arithmetic, are intermediate between $\mathsf{ACA}_0$
and $\mathsf{ATR}_0$.   There is much to be said for the IBSSM model over its finite version: we have remarked that the infinite version calculates power series functions, such as $\sin$, $e^{x}$. With  a little work one sees also that if any  differentiable function $f:\re \longrightarrow \re$ is IBSSM computable, then so is its derivative $f'$. 

On the other hand the class of sets that ITTM's compute form a Spector
class, and so we can bring to bear general results about such classes on the
ITTM semi-decidable, and eventually semi-decidable classes; their strength we saw was very strong: between $\mathsf{\Pi^1_2}$-$\mathsf{CA_{0}}$  and 
$\mathsf{\Pi^1_3}$-$\mathsf{CA_{0}}$. Finally the $\tmop{On}$-tape version of the ITTM,
gives us a new presentation of the constructible hierarchy as laid out by an
ordinary Turing program progressing throughout $\tmop{On}$ time.

\small

\bibliographystyle{plain}
\bibliography{settheory10h}

\begin{thebibliography}{10}

\bibitem{BeBuFr12}
A.~Beckmann, S.~Buss, and S-D. Friedman.
\newblock Safe recursive set functions.
\newblock {\em Centre de Recerca Matematica Document Series, Barcelona}, 2012.

\bibitem{BellCook92}
S.~Bellantoni and S.~Cook.
\newblock A new recursion-theoretic characterization of the poly-time
  functions.
\newblock {\em Computational Complexity}, 2:97--110, 1992.

\bibitem{BlShSm89}
L.~Blum, M.~Shub, and S.~Smale.
\newblock On a theory of computation and complexity over the real numbers.
\newblock {\em Notices of the American Mathematics Society (N.S.)},
  21(1):1--46, 1989.

\bibitem{B}
J.P. Burgess.
\newblock The truth is never simple.
\newblock {\em Journal of Symbolic Logic}, 51(3):663--681, 1986.

\bibitem{CFKMNW}
M.~Carl, T.~Fischbach, P.~Koepke, R.~Miller, M.~Nasfi, and G.~Weckbecker.
\newblock The basic theory of infinite time register machines.
\newblock {\em Archive for Mathematical Logic}, 49(2):249--273, 2010.

\bibitem{Da09}
B.~Dawson.
\newblock {\em Ordinal time {Turing} computation}.
\newblock PhD thesis, Bristol, 2009.

\bibitem{Fi03}
H.~Field.
\newblock A revenge-immune solution to the semantic paradoxes.
\newblock {\em Journal of Philosophical Logic}, 32(3):139--177, April 2003.

\bibitem{FrWe11}
S-D. Friedman and P.~D. Welch.
\newblock Hypermachines.
\newblock {\em Journal of Symbolic Logic}, 76(2):620--636, June 2011.

\bibitem{FrWe0?}
S-D. Friedman and P.D. Welch.
\newblock Two observations concerning infinite time {Turing} machines.
\newblock In I.~Dimitriou, editor, {\em BIWOC 2007 Report}, pages 44--47, Bonn,
  January 2007. Hausdorff Centre for Mathematics.
\newblock Also at http://www.logic.univie.ac.at/sdf/papers/joint.philip.ps.

\bibitem{Go65}
E.~Gold.
\newblock Limiting recursion.
\newblock {\em Journal of Symbolic Logic}, 30(1):28--48, Mar 1965.

\bibitem{HL}
J.D. Hamkins and A.~Lewis.
\newblock Infinite time {Turing} machines.
\newblock {\em Journal of Symbolic Logic}, 65(2):567--604, 2000.

\bibitem{HL2}
J.D. Hamkins and A.~Lewis.
\newblock {Post's} problem for supertasks has both positive and negative
  solutions.
\newblock {\em Archive for Mathematical Logic}, 41:507--523, 2002.

\bibitem{HaMi09}
J.D. Hamkins and R.~Miller.
\newblock Post's problem for ordinal register machines: an explicit approach.
\newblock {\em Annals of Pure and Applied Logic}, 160(3):302--309, September
  2009.

\bibitem{HaSe}
J.D. Hamkins and D.~Seabold.
\newblock Infinite time {Turing} machines with only one tape.
\newblock {\em Mathematical Logic Quarterly}, 47(2):271--287, 2001.

\bibitem{H}
H.G. Herzberger.
\newblock Notes on naive semantics.
\newblock {\em Journal of Philosophical Logic}, 11:61--102, 1982.

\bibitem{Kl59}
S.~C. Kleene.
\newblock Recursive quantifiers and functionals of finite type {I}.
\newblock {\em Transactions of the American Mathematical Society}, 91:1--52,
  1959.

\bibitem{Kl62b}
S.~C. Kleene.
\newblock Turing-machine computable functionals of finite type {I}.
\newblock In {\em Proceedings 1960 Conference on Logic, Methodology and
  Philosopy of Science}, pages 38--45. Stanford University Press, 1962.

\bibitem{Kl62a}
S.~C. Kleene.
\newblock Turing-machine computable functionals of finite type {II}.
\newblock {\em Proceedings of the London Mathematical Society}, 12:245--258,
  1962.

\bibitem{Kl63}
S.~C. Kleene.
\newblock Recursive quantifiers and functionals of finite type {II}.
\newblock {\em Transactions of the American Mathematical Society},
  108:106--142, 1963.

\bibitem{Kle07}
A.~Klev.
\newblock {\em Magister thesis}.
\newblock ILLC Amsterdam, 2007.

\bibitem{K05}
P.~Koepke.
\newblock Turing computation on ordinals.
\newblock {\em Bulletin of Symbolic Logic}, 11:377--397, 2005.

\bibitem{KK06}
P.~Koepke and M.~Koerwien.
\newblock Ordinal computations.
\newblock {\em Mathematical Structures in Computer Science}, 16.5:867--884,
  October 2006.

\bibitem{KoMi08}
P.~Koepke and R.~Miller.
\newblock An enhanced theory of infinite time register machines.
\newblock In A.~Beckmann {\em et al.}, editor, {\em Logic and the Theory of
  Algorithms}, volume 5028 of {\em Springer Lecture Notes Computer Science},
  pages 306--315. Swansea, Springer, 2008.

\bibitem{KoSe09}
P.~Koepke and B.~Seyfferth.
\newblock Ordinal machines and admissible recursion theory.
\newblock {\em Annals of Pure and Applied Logic}, 160(3):310--318, 2009.

\bibitem{KoSi06}
P.~Koepke and R.~Siders.
\newblock Computing the recursive truth predicate on ordinal register machines.
\newblock In A.~Beckmann {\em et al.}, editor, {\em Logical Approaches to
  Computational Barriers}, Computer Science Report Series, page~21. Swansea,
  2006.

\bibitem{KoWe11}
P.~Koepke and P.D. Welch.
\newblock A generalised dynamical system, infinite time register machines, and
  {$\Pi^1_1$-$\tmop{CA}_0$}.
\newblock In B.~L\, editor, {\em Proceedings of CiE 2011, Sofia}, 2011.

\bibitem{L}
B.~L{\"{o}}we.
\newblock Revision sequences and computers with an infinite amount of time.
\newblock {\em Journal of Logic and Computation}, 11:25--40, 2001.

\bibitem{Min67}
M.~Minsky.
\newblock {\em Computation: Finite and Infinite Machines}.
\newblock Prentice-Hall, 1967.

\bibitem{Pu65}
H.~Putnam.
\newblock Trial and error predicates and the solution to a problem of
  {Mostowski}.
\newblock {\em Journal of Symbolic Logic}, 30:49--57, 1965.

\bibitem{Rog}
H.~Rogers.
\newblock {\em Recursive Function Theory}.
\newblock Higher Mathematics. McGraw, 1967.

\bibitem{Sa66}
G.E. Sacks.
\newblock Post's problem, admissible ordinals and regularity.
\newblock {\em Transactions of the American Mathematical Society}, 124:1--23,
  1966.

\bibitem{Sa90}
G.E. Sacks.
\newblock {\em Higher Recursion Theory}.
\newblock Perspectives in Mathematical Logic. Springer Verlag, 1990.

\bibitem{SheStu63}
J.~Shepherdson and H.~Sturgis.
\newblock Computability of recursive functionals.
\newblock {\em Journal of the Association of Computing Machinery}, 10:217--255,
  1963.

\bibitem{Sho75}
R.~A. Shore.
\newblock Splitting an $\alpha$ recursively enumerable set.
\newblock {\em Transactions of the American Mathematical Society}, 204:65--78,
  1975.

\bibitem{Sho76}
R.~A. Shore.
\newblock The recursively enumerable $\alpha$-degrees are dense.
\newblock {\em Annals of Mathematical Logic}, 9:123--155, 1976.

\bibitem{Th54}
J.~Thomson.
\newblock Tasks and supertasks.
\newblock {\em Analysis}, 15(1):1--13, 1954-55.

\bibitem{W5}
P.D. Welch.
\newblock Minimality arguments in the infinite time {Turing} degrees.
\newblock In S.B.Cooper and J.K.Truss, editors, {\em Sets and Proofs: Proc.
  Logic Colloquium 1997, Leeds}, volume 258 of {\em London Mathematical Society
  Lecture Notes in Mathematics}. C.U.P., 1999.

\bibitem{W}
P.D. Welch.
\newblock Eventually {Infinite} {Time} {Turing} degrees: infinite time
  decidable reals.
\newblock {\em Journal of Symbolic Logic}, 65(3):1193--1203, 2000.

\bibitem{W8}
P.D. Welch.
\newblock Post's and other problems in higher type supertasks.
\newblock In B.~L{\"{o}}we, B.~Piwinger, and T.~R{\"{a}}sch, editors, {\em
  Classical and New Paradigms of Computation and their Complexity hierarchies,
  Papers of the Conference Foundations of the Formal Sciences III}, volume~23
  of {\em Trends in logic}, pages 223--237. Kluwer, Oct 2004.

\bibitem{W16}
P.D. Welch.
\newblock Ultimate truth {\em vis \`{a} vis} stable truth.
\newblock {\em Review of Symbolic Logic}, 1(1):126--142, June 2008.

\bibitem{W09}
P.D. Welch.
\newblock Characteristics of discrete transfinite {Turing} machine models:
  halting times, stabilization times, and normal form theorems.
\newblock {\em Theoretical Computer Science}, 410:426--442, January 2009.

\end{thebibliography}
\end{document}